\documentclass[onecolumn]{elsart3p}

\usepackage{amssymb}
\usepackage[utf8]{inputenc}
\usepackage[french,english]{babel}
\usepackage{times,amsmath,amsfonts,graphicx,scalefnt,fancyvrb}
\usepackage{color,listings} 
\newcommand{\vect}[1]{\mathbf{#1}}

\begin{document}

\begin{frontmatter}
\title{Fast computation of gradient and sentitivity in 13C metabolic flux analysis instationary experiments using the adjoint method}

\author{Stéphane Mottelet}
\address{Laboratoire de Math\'ematiques Appliqu\'ees de Compi\`egne, 
D\'epartement de G\'enie Informatique, 
Universit\'e de Technologie de Compi\`egne, 
BP 20529, 60205 COMPIEGNE CEDEX, FRANCE}
\ead{stephane.mottelet@utc.fr}


\author{}

\address{}

\begin{abstract}
\end{abstract}

\end{frontmatter}

\section{Motivation}
The overall dynamics of a CLE can be described by a cascade of differential equations of the 
following form (see e.g. \cite{wiechert}) :

\begin{equation}
\vect{X}_k(\vect{m})\vect{\dot x}_k = \vect{f}_k(\vect{v},\vect{x}_k,\dots,\vect{x}_1,\vect{x}_k^{input}),\;k=1\dots n,\;t\in[0,T]
\end{equation}

where the states $\vect{x}_k$ are functions of time $t$ and take their values in $\mathbb{R}^{n_k}$ (they represent the cumomer fractions of each metabolite) and the constant 
vectors $\vect{x}_k^{input}$ are vectors of  $\mathbb{R}^{n^{input}_k}$ depend on the labeling of the input substrates. The vector $\vect{v}\in\mathbb{R}^m$ denotes 
the unknown fluxes and $\vect{X}_k$  are diagonal matrices containing unknown pool sizes corresponding to to cumomer fractions of weight $k$.

The particular form taken by the functions $\vect{f}_k$ depends on the transition pattern of carbon atoms occuring for each reaction in the metabolic network. Writing down by hand the expression of these functions is quite easy for a small sample network but
becomes untractable for a realistic network. As far as numerical computations are concerned (direct problem solving or  identification) the real concern is to write some specific computer code computing these functions and their exact derivatives with respect of states
$\vect{x}_1,\dots,\vect{x}_n$ and $\vect{v}$, in a target language. The formal expression of the overall system could be interesting for testing the identifiability of the flux vector $\vect{v}$ and the pool sizes, but previous work shows that the size of realistic networks prevents the use of classical algorithms based on symbolic computations. 

Nowadays, the most efficient way of describing a metabolic network is to use the Systems Biology Markup Language (SBML, see \cite{sbml}, \cite{sbml2}), as it has become the de facto standard, used by a growing number of commercial or open source applications. The SBML markup language is a dedicated dialect of XML (see e.g. \cite{xml}) with a specific structure which allows to describe the different compartments, species, and the kinetics of reactions occuring between these species. Transformations can be applied to the SBML file describing the network, described in another XML dialect, the eXtendted Stylesheet Language (XSL), and the kind of transformations we are interested in, are those which will allow to generate the specific numerical code we need to solve the identification problem (stationnary and instationnary). The generated computer code is specific for each particular
metabolic network and associated CLE, an thus more efficient, readable and resusable, than a general application written to cover all possible cases. 

The target language which has been chosen is Scilab (see \cite{scilabbook,scilabweb})
because it is an open-source Matlab compatible language, allowing high-level programming together with compiled libraries with performant differential equation solvers, optimization routines and efficient sparse matrix algebra. The GUI of the final application is also described in XML, using another specific dialect called XMLlab (see \cite{xmllab}), which is available under the form of an official Scilab Toolbox (see Scilab www site). The description of the GUI is obtained with another pass of XSL transformations on the original SBML file describing the network. This GUI allows the biologist to enter the input data (label input of the substrate, known fluxes or a priori relationships betwen them, label observation), and lauch the optimization process solving the identification problem.

\begin{figure}[hb]
\begin{tabular}{ccc}
\begin{minipage}[c]{6cm}
$$
\includegraphics[width=4cm]{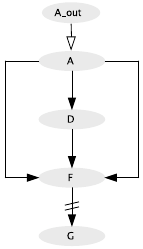}
$$
\end{minipage}
&\begin{minipage}[t]{3cm}
\begin{tabular}{rrcl}
$v_1$ : & A & $\to$ & F\\
 & \#ij & $\to$ & \#ij\\
$v_2$ : &A & $\to$ & D + D\\
  & \#ij & $\to$ &\#i + \#j\\
$v_3$ : &A & $\to$ & F\\
  & \#ij & $\to$ &\#ji
\end{tabular}
\end{minipage}

&\begin{minipage}[t]{3cm}
\begin{tabular}{rrcl}
$v_4$ : &D + D & $\to$ & F\\
  & \#i + \#j & $\to$ &\#ij\\
$v_5$ : &F & $\to$ & G\\
  & \#ij  & $\to$ &\#ij\\
$v_6$ : &A\_out& $\to$ & A\\
  & \#ij & $\to$ &\#ij
\end{tabular}
\end{minipage}
\end{tabular}
\caption{The branching network, with associated fluxes ${v}_1$, ${v}_2$, ${v}_3$, ${v}_4$, ${v}_5$, ${v}_6$ and the carbon atoms transition.}
\label{branching}
\end{figure}

Throughout this paper, we will use a very small example to illustrate our approach. This is the branching network used by Wiechert and Isermann
in \cite{wiechert2} (see Figure \ref{branching}).
The SBML file corresponding to this network can be found in Figures \ref{branching:listing1} and \ref{branching:listing2} in
the appendix section. As it can be seen, the description is very verbose. The added material describing some informations
on the CLE (label input, label observation, carbon atom mapping) are entered in the species and reaction notes, directly from 
the CellDesigner interface. In the future, we plan to develop a plugin to directly enter these information from within CellDesigner
and create XML annotation in the SBML file.
.
\section{Mathematical modelling in the stationnary case}

It has been shown in \cite{wiechert} that the actual state equation is in fact a succession of
linear ordinary differential equation, where for each $k$ the non homogeneous part $\vect{b}_k$ of the right handside depends on
$\vect{x}_1,\dots,\vect{x}_{k-1}$, giving the following cascade

\begin{equation}
\vect{X}_k(\vect{m})\vect{\dot x}_k(t) = \vect{M}_k(\vect{v})\vect{x}_k(t) +\vect{b}_k(\vect{v},\vect{x}_{k-1}(t),\dots,\vect{x}_1(t),\vect{x}_k^{input}),\;k=1\dots n,\;t>0.
\label{state2}
\end{equation}

Each component of vectors $\vect{x}_k(t)$ represents a cumomer fraction of weight $k$ of a 
given species (for a proper definition of
cumomer and cumomer weight see \cite{wiechert}). The constant vectors $\vect{x}_k^{input}$ 
contain the cumomer fractions of weight $k$ of species
which are input metabolites. The diagonal matrices $\vect{X}_k(\vect{m})$ depend on the stationnary 
concentrations of metabolites. The matrix $\vect{M}_k(\vect{v})$ and the vector $\vect{b}_k$ are constructed by
considering the balance equation for each cumomer of the vector $\vect{x}_k$. Constructing
these matrix by hand is a very tedious task but can be
automated if adequate data structures are used, in order to represent the metabolic network and the carbon transition map for each reaction (we
will explain later  how we deal with this particular information). 

\bigskip
In the stationnary case, the CLE is considering the asymptotic behaviour of the system (\ref{state2}), 

\begin{equation}
0 = \vect{M}_k(\vect{v})\vect{x}_k +\vect{b}_k(\vect{v},\vect{x}_{k-1}\dots,\vect{x}_1(t),\vect{x}_k^{input}),\;k=1\dots n,
\label{state:stat}
\end{equation}

in this case the states $\vect{x}_k$ do not depend on time anymore and the identification problem is restricted to
the determination of the flux vector $\vect{v}$ such that some cost function is minimized. For a given flux vector $\vect{v}$  this cost function can be classically defined as the squared norm of the difference betwen an observation vector $\vect{y}^{meas}\in \mathbb{R}^{n_{meas}}$ and the corresponding synthetic observation $\vect{y}(\vect{v})$ computed by solving the state equation (\ref{state:stat}) for the given value of the fluxes. 

In the following, we will consider that this observation is composed
of isotopomer and cumomer fractions of given species, which can always be computed as linear combination of cumomers, i.e. there exists $n$ non zero matrices $\vect{C}_1,\vect{C}_2,\dots,\vect{C}_n$ such that 
$$
\vect{y}(\vect{v})=\sum_{k=1}^n \vect{C}_k \vect{x}_k(\vect{v}),
$$
where we have used $\vect{x}_k(\vect{v})$ to denote the solutions of the state equation (\ref{state:stat}) for a given flux vector $\vect{v}$. 

We have also natural constraints on the flux vector which result of the particular structure of the metabolic network (the stoichiometric balances) and some other kinds of linear constraints on the fluxes, which can express that some of the fluxes are fixed, for example the flux
of input substrate, or some more specific information, e.g. some linear combination of fluxes which should be zero. Thus, the constraints on $\vect{v}$
take an affine form
$$
\vect{A}\vect{v}=\vect{w}.
$$
There is usually some measured extracellular fluxes $\vect{v}^{meas}$, which has to be compared with the actual value of these fluxes, which can be always be expressed as linear function of $\vect{v}$. The mostly used cost function is the Chi-Square function
$$
J(\vect{v})=\frac{1}{2}\left\Vert \vect{\sigma}^{-1}\left(\vect{y}(\vect{v})-\vect{y}^{meas}\right)\right\Vert^2 
+ \frac{1}{2}\left\Vert\vect{\alpha}^{-1}\left(\vect{E}\vect{v}-\vect{v}^{meas}\right)\right\Vert^2
$$
where $\vect{\sigma}$ and $\vect{\alpha}$ are diagonal positive definite matrices containing the standard deviation for each observation. Minimizing the Chi-Square function, under the hypothesis of gaussian distribution, 
is equivalent to maximizing the likelihood of measurements. 
\begin{equation}
\left\{
\begin{array}{rcl}
\vect{\vect{\hat v}}&=&\operatorname{arg}\min_{\vect{v}\in\mathbb{R}^m} J_{\varepsilon}(\vect{v}),\\
\vect{A}\vect{v}&=&\vect{w},\\
\vect{v} &\geq &0,
\end{array}
\right.
\tag{$P_\varepsilon$}
\end{equation}
where $\vect{y}(\vect{v})=\sum_{k=1}^n \vect{C}_k \vect{x}_k(\vect{v})$ and $\vect{x}_k(\vect{v})$, $k=1\dots n$ are the solutions of the
state equation (\ref{state:stat}). 

\subsection{Parametrisation of the admissible fluxes subspace}

The subspace of admissible fluxes is determined by the system of equations and inequations

\begin{equation}
\left\{
\begin{array}{rcl}
\vect{A}\vect{v}&=&\vect{w},\\
\vect{v} &\geq &0,
\end{array}
\right.
\label{sys:constr}
\end{equation}
In order to detect any redundancy or incompatibilities due to the eventual complimentary constraints
added by the user, an admissibility test is done on the system. We do it by solving a trivial linear
program, which allows to test if $w$ is in the range of $A$ and then if the subspace (\ref{sys:constr}) is
non-void.

Then there are two possibilities to obtain a parametrization : by computing an orthonormal basis $\{\vect{V}^1,\dots \vect{V}^r\}$ of the kernel of $\vect{A}$ (where 
$p=\operatorname{rank} \vect{A}$) and the minimum norm solution $\vect{v}_0$ of $\vect{A}\vect{v}=\vect{w}$, any $\vect{v}$
satisfying (\ref{sys:constr}) can be expressed as
$$
\vect{v=Vq+v}_0,
$$
where $\vect{q}$ is a vector of size $m-p$.

\bigskip
The classical parametrization, using the free fluxes, can be found  by computing the $\vect{Q}\vect{R}$ factorization of $\vect{A}$. There exists an $m\times m$ permutation matrix $\vect{P}=[\vect{P}_1,\vect{P}_2]$ and an orthogonal square matrix $\vect{Q}=[\vect{Q}_1,\vect{Q}_2]$, where $\vect{P}_1$ and $\vect{Q}_1$ represent the first $p$ columns of $\vect{P}$ and $\vect{Q}$, and a full rank $p\times p$ upper-triangular matrix $\vect{R}_1$ such that
$$
\vect{A}\vect{P}=\vect{Q}\left[\begin{array}{c|c}\vect{R}_1 & \vect{R}_2\\\hline 0 & 0\end{array}\right],
$$
where the lower right zero block is absent if $\vect{A}$ has full rank. The free fluxes in the vector are given by $\vect{q}=\vect{P}_2^\top \vect{v}$ and the complimentary dependent fluxes are given by $\vect{P}_1^\top \vect{v}$. Straightforward computations give the parametrization $\vect{v=Vq+v}_0$, where $\vect{q}$ has $m-p$ components and
$$
\vect{V}=\vect{P}_2-\vect{P}_1 \vect{R}_1^{-1}\vect{R}_2,~\vect{v}_0=\vect{P}_1 \vect{R}_1 ^{-1} \vect{Q}_1 ^\top \vect{w}.
$$

\bigskip Such parametrizations remove any redundancy in the constraints, and allow to identify which constraints in $\vect{V}\vect{q}+\vect{v}_0\geq 0$ are equality constraints blocking the value of some fluxes. Some of them are not specified in the initial system (\ref{sys:constr}) but are added {\em de facto} in the parametrization because of the implicit fluxes balancing constraints. This situation can be detected when a row of $\vect{V}_i$ is equal to zero. Hence, if we define the sets $I=\{i,\;\vect{V}_i\neq 0\}$ and the set $D$ containing the indices of dependent fluxes, the parametrized optimization problem takes the form
\begin{equation}
\left\{
\begin{array}{rcl}
\vect{\hat q}&=&\operatorname{arg}\min_{\vect{q}\in\mathbb{R}^{m-p}} J_{\varepsilon}(\vect{V}\vect{q}+\vect{v}_0),\\
q_i&\geq& 0,\;i=1\dots m-p,\\
\vect{V}_i\vect{q}+(\vect{v}_0)_i&\geq& 0,\;i\in I\cap D,
\end{array}
\right.
\label{pparam}
\end{equation}
and we have $\vect{\vect{\hat v}}=\vect{V}\vect{\hat q}+\vect{v}_0$. The gradient of the parametrized cost function can be expressed via the chain rule as
$$
\left(\frac{d}{d\vect{q}}J_{\varepsilon}(\vect{V}\vect{q}+\vect{v}_0)\right)^\top=\vect{V}^\top\nabla J_\varepsilon(\vect{V}\vect{q}+\vect{v}_0).
$$
The type of parametrization (orthogonal or free fluxes)
used in the optimization does not seem to influence the conditioning of the algorithm, so the free fluxes parametrization is used, because of its biological interpretation.

\subsection{Identifiability and regularization}
When $\varepsilon=0$ the constraints on $\vect{q}$ are not enough to  ensure existence of a solution
because the problem may be unbounded. In fact, existence and unicity of a solution will occur if the fluxes are identifiable. A general discussion about this subject can be found in \cite{wiechert2}, where the authors propose an algorithm based on integer arithmetics to test the structural identifiability of metabolic networks. The most encoutered problematic situation corresponds to bidirectional reactions, such as 
$$
v_1:\mathrm{A\to B},~v_2:\mathrm{B\to A}
$$
where $v_1-v_2$ (the net flux) is identifiable but $v_1$ and $v_2$ are not individualy identifiable. The  counterpart of such a situation is that the optimal $v_1$ and $v_2$ tend to infinity when $\varepsilon\to 0$, and the cost function $J_{\varepsilon}$ is ill-conditioned when $\varepsilon$ is too small, leading to convergence problems in the optimization phase. The change of variables proposed in \cite{wiechertnonstat} considers the net flux $v_{net}=v_1-v_2$ and the exchange fluxes $v_{xch}=\min(v_1,v_2)$ and a \og compacification\fg{} of $v_{xch}$ defined by 
$$
v_{[0,1]}=\frac{v_{xch}}{\beta+v_{xch}},
$$ 
where $\beta>0$. The above change of variables maps $[0,+\infty[$ to $[0,1[$ and thus is interesting from a numerical point of view. Although these new variables make sense from a metabolic point of view, it remains that the overall mapping from $(v_1,v_2)$ to $(v_{net},v_{[0,1]})$ is not differentiable and thus needs to be approximated. A more systematic approach is proposed in \cite{yang} where all free fluxes $q_i\geq 0$ are mapped to $r_i\in [0,1[$ with the change of variables $\vect{q}=\vect{q}(\vect{r})$, where
$$
q_i=\beta\frac{r_i}{1-r_i},\;i=1\dots m-p,
$$
where $\beta>0$ is a scaling parameter. In this case, the inequality constraints in (\ref{pparam}) become non linear and the new optimization problem is
\begin{equation}
\left\{
\begin{array}{rcl}
\vect{\hat r}&=&\operatorname{arg}\min_{\vect{r}\in\mathbb{R}^{m-p}} J_{\varepsilon}(\vect{V}\vect{q(r)}+\vect{v}_0),\\
1-\delta\geq r_i&\geq& 0,\;i=1\dots m-p,\\
\vect{V}_i\vect{q(r)}+(\vect{v}_0)_i&\geq& 0,\;i\in I\cap D,
\end{array}
\right.
\label{pparamcompact}
\end{equation}
where $\delta>0$ can be arbitrary small. The gradient of the cost function is given by
$$
\left(\frac{d}{d\vect{r}}J_{\varepsilon}(\vect{V}\vect{q(r)}+\vect{v}_0)\right)^\top=\vect{q'(r)}^{\top}\vect{V}^\top\nabla J_\varepsilon(\vect{V}\vect{q(r)}+\vect{v}_0).$$
\subsection{Multiple experiences}

In the following we will also consider the case where multiple CLE are done with the same metabolic network but with diffferent
labeling of the input metabolites, given by $\vect{x}^{input,i}$ for $i=1\dots n_{exp}$. Thus, we will consider the cost function

$$
J_{\varepsilon}(\vect{v})=\frac{1}{2}\sum_{i=1}^{n_{exp}}\left(\left\Vert
\vect{\sigma}^{-1}\left(\vect{y}(\vect{v},\vect{x}^{input,i})-\vect{y}^{meas,i}\right)\right\Vert^2 + \left\Vert
\vect{\alpha}^{-1}\left(\vect{E}\vect{v}-\vect{v}_{obs,i}\right)\right\Vert^2\right) + \frac{\varepsilon}{2} \Vert\vect{v}\Vert^2
$$
where $\vect{y}^{meas,i}$ is the observation of labeled material for experience $i$, $\vect{v}_{obs,i}$ is the vector of measured extracellular fluxes and
$$\vect{y}(\vect{v},\vect{x}^{input,i})=\sum_{k=1}^n \vect{C}_k \vect{x}_k(\vect{v},\vect{x}^{input,i})$$ and $\vect{x}_k(\vect{v},\vect{x}^{input,i})$, $k=1\dots n$ is the solution of
\begin{equation}
0 = \vect{M}_k(\vect{v})\vect{x}_k +\vect{b}_k(\vect{v},\vect{x}_{k-1},\dots,\vect{x}_1(t),\vect{x}_k^{input,i}),\;k=1\dots n,
\label{state:stat:mult}
\end{equation}

\subsection{Computation of the gradient of the cost function}

The computation of the gradient of $J(\vect{v})$ needs the derivative of $\vect{x}(\vect{v})$ with respect to $\vect{v}$. In the stationnary case, it makes sense to compute this derivative by implicit differentiation of the state equation (\ref{state:stat}). To this purpose,
we adopt the notation 
\begin{equation}
\vect{f}_k(\vect{v},\vect{x}_k,\dots,\vect{x}_1,\vect{x}_k^{input,i})=\vect{M}_k(\vect{v})\vect{x}_k +\vect{b}_k(\vect{v},\vect{x}_{k-1},\dots,\vect{x}_1,\vect{x}_k^{input,i}) and
\label{notation:fk}
\end{equation}
we denote by $\vect{x}^i(\vect{v})$ the solution of 
$$
\vect{f}_k(\vect{v},\vect{x}_k,\dots,\vect{x}_1,\vect{x}_k^{input,i})=0,\;k=1\dots n.
$$
By differentiating these equations with respect to $\vect{v}$, when $\vect{x}=\vect{x}^i(\vect{v})$, we obtain for $k=1\dots n$
\begin{equation}
0=\frac{d\vect{f}_k}{d\vect{v}}(\vect{v},\vect{x}_k,\dots,\vect{x}_1,\vect{x}_k^{input,i})=
\frac{\partial \vect{f}_k}{\partial\vect{v}}(\vect{v},\vect{x}_k,\dots,\vect{x}_1,\vect{x}_k^{input,i})
+\sum_{l=1}^{k}\frac{\partial \vect{f}_l}{\partial \vect{x}_l}(\vect{v},\vect{x}_k,\dots,\vect{x}_1,\vect{x}_k^{input,i}) \frac{\partial \vect{x}^i_l}{\partial\vect{v}}.
\end{equation}
Since $\vect{f}_k$ is linear with respect to $\vect{x}_k$ for fixed $\vect{v}$, we can determine
$\frac{\partial \vect{x}_k}{\partial\vect{v}}$ as the solution of a linear system of equations, whose right hand side is a function
of $\vect{x}_l$ and $\frac{\partial \vect{x}_l}{\partial\vect{v}}$ for $l=1\dots k$:
\begin{equation}
\vect{M}_k(\vect{v})\frac{\partial \vect{x}^i_k}{\partial\vect{v}}=\frac{\partial \vect{f}_k}{\partial\vect{v}}(\vect{v},\vect{x}_k,\dots,\vect{x}_1,\vect{x}_k^{input,i})
-
\sum_{l=1}^{k-1}\frac{\partial \vect{b}_k}{\partial \vect{x}_l}(\vect{v},\vect{x}_k,\dots,\vect{x}_1,\vect{x}_k^{input}) \frac{\partial \vect{x}^i_l}{\partial\vect{v}},\;k=1\dots n.
\label{cascade:deriv}
\end{equation}
Hence, the key ingredients in the computation of the derivatives of $\vect{x}_k$ are the derivatives $\frac{\partial \vect{b}_k}{\partial \vect{x}_l}$
for $l<k$ and the derivatives $\frac{\partial \vect{f}_k}{\partial\vect{v}}$. This will be one of the main tasks of the automatically generated
computer code, together with the assembly of matrices $\vect{M}_k(\vect{v})$. Since the gradient of the cost function $J(\vect{v})$ will be required at each
iteration of the optimization algorithm, these matrix will be assembled as sparse matrices in order to speed up the computations,
particularly the resolution of the linear systems (\ref{cascade:deriv}).

\bigskip
The final computation of the derivative gives
$$
\frac{d J_\varepsilon(\vect{v})}{d\vect{v}}=\left(\vect{E}\vect{v}-\vect{v}_{obs,i}\right)^\top \vect{\alpha}^{-2}\vect{E}+\varepsilon \vect{v}^\top+\sum_{i=1}^{n_{exp}}\left(\vect{y}(\vect{v},\vect{x}^{input,i})-\vect{y}^{meas}\right)^\top\vect{\sigma}^{-2}\sum_{k=1}^n\vect{C}_k\frac{\partial
\vect{x}^i_k}{\partial\vect{v}}.
$$

\section{Architecture of the computer code generation algorithms}

The most innovative aspect of this work is the choice of the techniques to generate the code : from the
original SBML file edited under Cell Designer (or any other SBML compliant software), only XSL (eXetended Stylesheet Language)
transformations are used to generate the Scilab code computing the specific objects for a given Carbon Labeling Experiment. The way
transformations are done is described in XSL stylesheets, written in anthor XML dialect. XSL is very different from the 
typical programming languages in use today. One question that's being asked frequently is : what kind of programming language is actually
XSLT ? Until now, the authoritative answer from some of the best specialists is that XSLT is a declarative (as opposed to imperative) language.
The XSL stylesheets are thus very explicit, human readable, and easy to debug and maintain.

\bigskip In the whole process which maps the original SBML file to the computer code and the graphical user interface, successive XSL transformations occur. The first set of transformations aims to translate all the 
specific information about the CLE into XML markup which can be later used, for example, the carbon atom mapping of each reaction (this step is described in Appendix A). Then two different main paths are followed : 

\bigskip \begin{enumerate}
\item The first series of transformations is dedicated to the computer code generation :
\begin{enumerate}
\smallskip\item  A main assembly loop is processed, which for each weight $k$, generates, for $j=1\dots n_j$, some intermediate XML markup declaring the contribution of each cumomer fraction $(\vect{x}_k)_j$ to the matrices
$$\vect{M}_k(\vect{v}),\,\vect{b}_k(\vect{v}),\frac{\partial \vect{f}_k}{\partial\vect{v}},\left(\frac{\partial\vect{b}_k}{\partial \vect{x}_l}\right)_{l=1\dots k-1}.$$
The contributions to $\vect{M}_k(\vect{v})$ are functions of the flux vector $\vect{v}$ only, but the contributions to other
matrices are functions of lower weight cumomer components and eventually of $\vect{x}_k^{input}$, respecting the weight preservation property (see \cite{wiechert}). This intermediate step is event-driven, i.e. contributions to the matrices are dumped in the order they occur. The contributions are gathered for each matrix in the subsequent transformation which produces the Scilab code.

\medskip\item The Scilab code solving the cascaded linear systems is generated. For each weight $k$, the matrix $\vect{M}_k(\vect{v})$ is stored as a sparse matrix and its sparse $LU$ factorization is computed (the sparse triangular factors are also retained because they are also needed to compute the derivatives). The linear system giving $\vect{x}_k$ is then solved by using the precomputed sparse $LU$ factors. Then we compute the matrices ${\partial \vect{f}_k}/{\partial\vect{v}}$ and $\left({\partial\vect{b}}/{\partial \vect{x}_i}\right)_{i=1\dots k-1}$, which need the previously computed cumomer vectors $\vect{x}_1,\dots,\vect{x}_{k-1}$ and finally
solve the linear system giving ${\partial \vect{x}_k}/{\partial\vect{v}}$.
\end{enumerate}
\item The second series of transformation aims to build the specific graphical user interface for the given metabolic network. To this purpose, an XML file conforming to the XMLlab DTD is generated. The structure of the interface is described in a high level way : there are given sections of the interface, each of these being dedicated to different purposes. The first section is hosting all the fluxes, the second section is hosting the fluxes observations with associated standard deviation, the third section hosts together the label measurements and the label output corresponding to the current fluxes. This is the place where the user can compare the original measurement and the reconstructed measurement after the identification process. The fourth section displays all components of the cumomer vector, sorted by weight and by species name. The last section is reserved to the parameters of the identification method (maximum iteration, regularization parameter, and so on). The structure of the original SBML file, enriched with the specific annotations in
the sbml namespace (see Appendix A), allows to perform this step in a very straightforward way.  
\end{enumerate}
The generated Scilab code computing the cumomers vector $\vect{x}$, the derivative matrices, the cost function and its gradient are given on Figures \ref{fig:branchingScilab1} and \ref{fig:branchingScilab2} in Appendix A.

\section{Numerical results}

\section{Mathematical modelling in the unstationnary case}

In the unstationnary case, the CLE is done when the system (\ref{state2}) has not reached
its asymptotic behaviour. The measurements
$$
\vect{y}^{meas,j},\;j=1\dots n_t,
$$
are done at different values $t_j$ of time, and we can make the hypothesis that the final time $T$ in (\ref{state2}) is equal to the final observation time i.e $T=t_{n_t}$. The fundamental difference with the stationnary case is that the stationnary concentration of metabolites is also an unknown, i.e. the vector $\vect{m}$ is also to be determined. 

The cost function takes the form
$$
J_\epsilon(\vect{v},\vect{m})=\frac{1}{2}\sum_{j=1}^{n_t} \left(\left\Vert \vect{\sigma}^{-1}\left(\vect{y}(t_j,\vect{v},\vect{m})-\vect{y}^{meas,j}\right)\right\Vert^2 \right)
+ \frac{1}{2}\left\Vert\vect{\alpha}^{-1}\left(\vect{E}\vect{v}-\vect{v}^{meas}\right)\right\Vert^2,
$$
where
$$\vect{y}(t_j,\vect{v},\vect{m})=\sum_{k=1}^n \vect{C}_k \vect{x}_k(t_i,\vect{v},\vect{m}),$$
and $\vect{x}_k(t,\vect{v},\vect{m})$, $k=1\dots n$ are the solutions of the state equation for a given pair $(\vect{v},\vect{m})$ of fluxes and pool sizes :

\begin{equation}
\vect{X}_k(\vect{m})\vect{\dot x}_k(t) = \vect{M}_k(\vect{v})\vect{x}_k(t) +\vect{b}_k(\vect{v},\vect{x}_{k-1}(t),\dots,\vect{x}_1(t),\vect{x}_k^{input}),\;k=1\dots n,\;t>0.
\label{cascadenonstat}
\end{equation}

The minimization problem is the the following : 

\begin{equation}
\left\{
\begin{array}{rcl}
(\vect{\hat v},\vect{\hat m})&=&\operatorname{arg}\min_{\vect{v},\vect{m}} J_{\varepsilon}(\vect{v},\vect{m}),\\
\vect{A}\vect{v}&=&\vect{w},\\
\vect{v} &\geq &0,\\
\vect{m} & \geq& 0.
\end{array}
\right.
\tag{$P^u_\varepsilon$}
\end{equation}

The main difficulty is the computation of the gradient of $J$, which can be done by using the sensitivity matrices $\frac{\partial \vect{x}_k(t)}{\partial\vect{v}}$, computed as the solution of a system of differential equations obtained by implicit differentiation of (\ref{cascadenonstat}) as in the stationnary case. The problem is that this approach is computationnaly intensive (see e.g. \cite{wiechertnonstat,Noh,Noh2006554}) because it involves a cascade of differential equations where the state is a matrix (instead of a vector).

A more suitable approach for non-stationnary problems is to use the adjoint state method (see \cite{plessix,chavent,lionsmagenes}). If the number of parameters of interest (the fluxes) exceeds the number of model outputs for which the sensitivity is desired, the adjoint method is more efficient than traditional direct methods of calculating sensitivities. The gradient of $J$ can be computed at the same cost as the state equation (\ref{cascadenonstat}). A far as the statistical evalution of identified fluxes is concerned, the sensitivity of the model output $\vect{y}(t,\vect{v},\vect{m})$ with respect to $\vect{v}$ can be obtained with a cost of $n_{meas}$ state equations. Hence, the adjoint method will always outperform the sensitivity method for the computation of the gradient. For the output sensitivity, the dimension of $\vect{y}(t,\vect{v},\vect{m})$ observations is generaly lower than the number of fluxes, so the same method should be used. The adjoint state method is best understood in continuous time and the next section is devoted to this presentation. Section \ref{sect:adjdiscr} will detail its practical implementation in discrete time. 

\subsection{The adjoint equation in continuous time}
In order to simplify the presentation of the results, we will adopt the block notation  $\vect{x}=\left(\vect{x}_1;\vect{x}_2;\dots \vect{x}_n\right)$ for the overall cumomer vector, and for the state equation we will write
\begin{equation}
\vect{X}(\vect{m})\vect{\dot x}(t)-\vect{f}(\vect{x}(t),\vect{v})=0,\;t\in[0,T[,
\label{state}
\end{equation}
where $\vect{f}(\vect{x},\vect{v})=\left(\vect{f}_1(\vect{x},\vect{v});\dots\vect{f}_n(\vect{x},\vect{v})\right)$ where $\vect{f}_k$ is defined by (\ref{notation:fk}).
We also define $\vect{C}=[\vect{C}_1,\vect{C}_2,\dots,\vect{C}_n]$ so that $\vect{y}(t)=\vect{C}\vect{x}(t)$. Without loss of generality, we consider only one measurement at final time $T$. 

\bigskip The adjoint state method allows to compute the total derivative with respect to $\vect{v}$ and $\vect{m}$ of a given function $I(\vect{x}(\vect{v},\vect{m}))\in\mathbb{R}$ where $\vect{x}(\vect{v},\vect{m})$ is the solution of the state equation (\ref{state}). If the gradient of $J_{\epsilon}$ is to be computed then we will take
\begin{equation}
I(\vect{x})=\frac{1}{2}\left\Vert \vect{\sigma}^{-1}\left(\vect{Cx}(T)-\vect{y}^{meas}\right)\right\Vert^2.
\label{I_for_J}
\end{equation}
In the following, we do not consider the quantities in $J_{\varepsilon}$ depending explicitely on $\vect{v}$ hence, we define $J(\vect{v},\vect{m})=I(\vect{x}(\vect{v},\vect{m}))$. Let us define the Lagrangian
\begin{equation}
\begin{split}
L(\vect{x},\vect{p},\vect{v},\vect{m})&=I(\vect{x})+\int_0^T \vect{p}(t)^\top\left(\vect{X}(\vect{m})\vect{\dot x}(t)-\vect{f}(\vect{x}(t),\vect{v})\right)dt,
\end{split}
\end{equation}
where the adjoint state $\vect{\vect{p}}=\left(\vect{p}_1;\vect{p}_2;\dots \vect{p}_n\right)$ has the same block structure as $\vect{x}$. The first remark that can be done is that when $\vect{x}$ is the solution of the state equation (\ref{state}), we have
$$
L(\vect{x}(\vect{v},\vect{m}),\vect{p},\vect{v},\vect{m})=J(\vect{v},\vect{m}),
$$
and when we express the total derivative of $J(\vect{v},\vect{m})$ e.g. with respect to $\vect{v}$ we have
\begin{equation}
\frac{dJ(\vect{v},\vect{m})}{d\vect{v}}=\frac{\partial L}{\partial \vect{x}}(\vect{x}(\vect{v},\vect{m}),\vect{p},\vect{v},\vect{m})\frac{\partial \vect{x}(\vect{v},\vect{m})}{\partial\vect{v}}+\frac{\partial L}{\partial\vect{v}}(\vect{x}(\vect{v},\vect{m}),\vect{p},\vect{v},\vect{m}).
\label{relat:adj}
\end{equation}
The idea of the adjoint state technique is to compute $\vect{p}$ such that $\frac{\partial L}{\partial \vect{x}}=0$. Then the remaining part of the derivative can be computed in a straightforward way. This adjoint equation is a (backward in time) differential equation given by
\begin{align}
\vect{X}(\vect{m})\vect{\dot p}(t)&=-\left(\frac{\partial{\vect{f}}}{\partial \vect{x}}(\vect{x}(t),\vect{v})\right)^\top \vect{p}(t),\;t\in[0,T[,
\label{adj}
\end{align}
with the final conditon
\begin{equation}
\vect{p}(T)=-\vect{C}^\top\vect{\sigma}^{-2}\left(\vect{C}\vect{x}(T)-\vect{y}^{meas}\right).
\label{finalcondadj}
\end{equation}
Because of the block triangular structure of $\frac{\partial{\vect{f}}}{\partial \vect{x}}$, the adjoint equation has also a cascade structure, but in reverse order, i.e. $\vect{p}_n$ is obtained at first and $\vect{p}_1$ at last :
\begin{align}
\label{finalcondk}\vect{p}_k(T)&=-\vect{C}_k^\top\vect{\sigma}^{-2}\left(\vect{C}\vect{x}(T)-\vect{y}^{meas}\right),\\
\vect{X}_k(\vect{m}) \vect{\vect{\dot p}}_k(t)&=-\vect{M}_k^\top\vect{p}_k(t)+\sum_{l=k+1}^{n}\left(\frac{\partial \vect{b}_l}{\partial \vect{x}_k}\right)^\top \vect{p}_l(t),\;t\in]0,T[,
\label{adjk}
\end{align}
for $k=1\dots n$.

\bigskip
When for a given pair $(\vect{v},\vect{m})$ the state equations and the adjoint state equations are solved, then the gradient of $J$ can be readily computed by using (\ref{relat:adj})
\begin{equation*}
\frac{dJ(\vect{v},\vect{m})}{d\vect{v}}=\frac{\partial L}{\partial\vect{v}}(\vect{x},\vect{p},\vect{v},\vect{m})=\sum_{k=1}^n\int_0^T \vect{p}_k(t)^\top\frac{\partial \vect{f}_k}{\partial\vect{v}}(\vect{v},\vect{x}(t))\,dt.
\end{equation*}
We will give the derivative with respect to $m$ in the next section.

\begin{rem}\rm
The output sensitivity $\frac{\partial \vect{y}(T)}{\partial\vect{v}}$ can be computed in the same way by taking $I(\vect{x})=\vect{C}\vect{x}(T)$. In this case, we have
$$
\frac{\partial \vect{y}(T)}{\partial\vect{v}}=\sum_{k=1}^n\int_0^T \vect{p}_k(t)^\top\frac{\partial \vect{f}_k}{\partial\vect{v}}(\vect{v},\vect{x}(t))\,dt,
$$
where the final condition (\ref{finalcondk}) is replaced by $\vect{p}_k(T)=-\vect{C}_k^\top$ and the adjoint equation (\ref{adjk}) is unchanged, but $\vect{p}_k(t)$ is a matrix of size $n_k\times n_{meas}$.
\end{rem}

\subsection{The adjoint equation in discrete time}
\label{sect:adjdiscr}
The previous sections shows that once the state equation is solved, the gradient of $J$ can be computed at the cost of one more differential equation (\ref{adj}), which has to be compared to the cost of computing the sensitivity function. But the practical implementation needs to reconsider this approach in discrete time, since we cannnot just approximate independently the continuous state and the continuous adjoint state equation, i.e. the discretized adjoint must be the adjoint of the discretize state. This is a reason why high order integration schemes are seldom used in adjoint codes written by hand (otherwise automatic differentiation can be used, see \cite{Bischof}) since the discrete adjoint is obtained by formal derivation of the state integration scheme. Since we have to consider that the state equation could be stiff because of eventual large values of fluxes, a good compromise is the implicit trapezoidal rule, which is of order 2. Hence, we consider a discretization of the interval $[0,T]$ by considering $t_i=(i-1)h$, for $i=1\dots N$ and $h=T/(N-1)$, and we denote by $\vect{x}^i$ the approximation of $\vect{x}(t_i)$. The implicit trapezoidal rule applied to equation (\ref{state}) gives
\begin{equation}
\vect{D(m)}(\vect{x}^{i+1}-\vect{x}^{i})-\frac{h}{2}(\vect{f}(\vect{x}^{i+1},\vect{v})+\vect{f}(\vect{x}^{i},\vect{v}))=0,\;i=1\dots N-1,
\label{state:discr}
\end{equation}
and we still denote by $\vect{x}(\vect{v},\vect{m})$ the solution of (\ref{state:discr}). We consider
that for each measurement time $\left\{\tau_j\right\}_{1\leq j\leq n_t}$ there exists $\theta(j)$ such that $\tau_j=t_{\theta(j)}$, with $\theta(n_t)=N$, and we define the cost function $J(\vect{v},\vect{m})=I(\vect{x}(\vect{v},\vect{m})$ where
$$I(\vect{x})=\frac{1}{2}\sum_{j=1}^{n_t}\left\Vert \vect{\sigma}^{-1}\left(\vect{C}\vect{x}^{\theta(j)}-\vect{y}^{meas,j}\right)\right\Vert^2.$$
The discrete Lagrangian is defined by
\begin{equation}
L(\vect{x},\vect{p},\vect{v},\vect{m})=I(\vect{x})+\sum_{i=1}^{N-1}(\vect{p}^i)^\top
(\vect{D(m)}(\vect{x}^{i+1}-\vect{x}^{i})-\frac{h}{2}(\vect{f}(\vect{x}^{i+1},\vect{v})+\vect{f}(\vect{x}^{i},\vect{v}))),
\end{equation}
where $\vect{p}^i$ is the adjoint state for time $i$. Straightforward computations show that the adjoint equation is given by
\begin{equation}
\left(\vect{X}(\vect{m})-\frac{h}{2}\frac{\partial \vect{f}}{\partial \vect{x}}(\vect{x}^i,\vect{v})^\top\right)\vect{p}^{i-1}=
\left(\vect{X}(\vect{m})+\frac{h}{2}\frac{\partial \vect{f}}{\partial \vect{x}}(\vect{x}^i,\vect{v})^\top\right)\vect{p}^i-\left(\frac{\partial I}{\partial \vect{x}^i}\right)^{\top},\;1<i<N,
\end{equation}
with the final condition
\begin{equation}
\left(\vect{X}(\vect{m})-\frac{h}{2}\frac{\partial \vect{f}}{\partial \vect{x}}(\vect{x}^N,\vect{v})^\top\right)\vect{p}^{N-1}=-\left(\frac{\partial I}{\partial \vect{x}^N}\right)^{\top},
\end{equation}
where $\frac{\partial I}{\partial \vect{x}^i}=0$ if $\theta(j)\neq i$ for all $j=1\dots n_t$ and
$$
\frac{\partial I}{\partial \vect{x}^i}=-\left(\vect{C}\vect{x}^i-\vect{y}^{meas,\theta^{-1}(i)}\right)^\top\sigma^{-2}\vect{C}
$$
otherwise. Once the state and the adjoint state equations are solved, the gradient is given by
$$
\left(\frac{dJ(\vect{v},\vect{m})}{d\vect{v}}\right)^\top=-\frac{h}{2}\sum_{i=1}^{N-1}\left(\frac{\partial \vect{f}}{\partial
\vect{v}}(\vect{x}^{i+1},\vect{v})+\frac{\partial \vect{f}}{\partial \vect{v}}(\vect{x}^{i},\vect{v})\right)^\top \vect{p}^i,
$$
and
$$
\left(\frac{dJ(\vect{v},\vect{m})}{d\vect{m}}\right)^\top=\sum_{i=1}^{N-1}
\left(\frac{\partial}{\partial\vect{m}}\left(\vect{X}(\vect{m})\left(\vect{x}^{i+1}-\vect{x}^{i}\right)\right)\right)^{\top}\vect{p}^i,
$$
where for a given vector $\vect{z}$ the matrix $\frac{\partial}{\partial \vect{m}}(\vect{X}(\vect{m})\vect{z})$ is defined by
\begin{equation}
\left(\frac{\partial}{\partial \vect{m}}(\vect{X}(\vect{m}) \vect{z})\right)_{ij}=\left\{\begin{array}{rl}
z_i,&\mbox{ if the cumomer fraction }z_i\mbox{ belongs to metabolite }j,\\
0,&\mbox{ otherwise.}
 \end{array}
 \right.
 \label{def:derivm}
\end{equation}

\appendix
\section{Workflow of SBML markup to Scilab code}

\begin{figure}[htb]
\begin{lstlisting}
<?xml version="1.0" encoding="UTF-8"?>
<sbml level="2" version="1" xmlns="http://www.sbml.org/sbml/level2"
      xmlns:ns="http://www.sbml.org/sbml/level2"
      xmlns:celldesigner="http://www.sbml.org/2001/ns/celldesigner">
 
  <model id="branching">
    <listOfCompartments>
      <compartment id="default" />
    </listOfCompartments>

    <listOfSpecies>
      <species compartment="default" id="A"/>
      <species compartment="default" id="D"/>
      <species compartment="default" id="F">
        <notes>
          <html xmlns="http://www.w3.org/1999/xhtml">
            <body> LABEL_MEASUREMENT 1x,x1,11 </body>
          </html>
        </notes>
      </species>

      <species compartment="default" id="G"/>

      <species compartment="default" id="A_out">
        <notes>
          <html xmlns="http://www.w3.org/1999/xhtml">
            <body> LABEL_INPUT 01,10,11 </body>
          </html>
        </notes>
      </species>
    </listOfSpecies>

    <listOfReactions>

      <reaction id="v1" reversible="false">
        <notes>
          <html xmlns="http://www.w3.org/1999/xhtml">
            <body> IJ &gt; IJ </body>
          </html>
        </notes>
        <listOfReactants>
          <speciesReference species="A"/>
        </listOfReactants>
        <listOfProducts>
          <speciesReference species="F"/>
        </listOfProducts>
      </reaction>

      <reaction id="v2" reversible="false">
        <notes>
          <html xmlns="http://www.w3.org/1999/xhtml">
            <body> IJ &gt; I+J </body>
          </html>
        </notes>
        <listOfReactants>
          <speciesReference species="A" />
        </listOfReactants>
        <listOfProducts>
          <speciesReference species="D" stoichiometry="2.0" />
        </listOfProducts>
      </reaction>
\end{lstlisting}
\caption{The XML file describing the branching network, with added information concerning the 
carbon mapping for each reaction, and the detail of input and measured label (lines 1 to 61)}
\label{branching:listing1}
\end{figure}
\begin{figure}[htb]
\begin{lstlisting}[firstnumber=62]
      <reaction id="v3" reversible="false">
        <notes>
          <html xmlns="http://www.w3.org/1999/xhtml">
            <body> IJ &gt; JI </body>
          </html>
        </notes>
        <listOfReactants>
          <speciesReference species="A"/>
        </listOfReactants>
        <listOfProducts>
          <speciesReference species="F"/>
        </listOfProducts>
      </reaction>

      <reaction id="v4" reversible="false">
        <notes>
          <html xmlns="http://www.w3.org/1999/xhtml">
            <body> I+J &gt; IJ </body>
          </html>
        </notes>
        <listOfReactants>
          <speciesReference species="D" stoichiometry="2"/>
        </listOfReactants>
        <listOfProducts>
          <speciesReference species="F"/>
        </listOfProducts>
      </reaction>

      <reaction id="v5"reversible="false">
        <notes>
          <html xmlns="http://www.w3.org/1999/xhtml">
            <body> IJ &gt; IJ </body>
          </html>
        </notes>
        <listOfReactants>
          <speciesReference species="F"/>
        </listOfReactants>
        <listOfProducts>
          <speciesReference species="G"/>
        </listOfProducts>
      </reaction>

      <reaction id="v6" reversible="false">
        <notes>
          <html xmlns="http://www.w3.org/1999/xhtml">
            <body> IJ &gt; IJ </body>
          </html>
        </notes>
        <listOfReactants>
          <speciesReference species="A_out"/>
        </listOfReactants>
        <listOfProducts>
          <speciesReference species="A" />
        </listOfProducts>
      </reaction>
	  
    </listOfReactions>
	
  </model>
</sbml>
\end{lstlisting}
\caption{The XML file describing the branching network, with added information concerning the 
carbon mapping for each reaction, and the detail of input and measured label (lines 62 to 121).}
\label{branching:listing2}
\end{figure}

The specific annotation in the private namespace \verb1xmlns:smtb="http://www.utc.fr/sysmetab"1 namespace concerning the carbon atom mapping of each reaction is done as follows : for example, for the reaction corresponding to flux $v_2$,
\begin{center}
\begin{minipage}[t]{3cm}
\begin{tabular}{rrcl}
$v_2$ : &A & $\to$ & D + D\\
  & \#ij & $\to$ &\#i + \#j
\end{tabular}
\end{minipage}
\end{center}
some specific markup is generated from the string \verb1IJ &gt; I+J1 found in the reaction \verb1<notes>1 element, as depicted in Figure \ref{fig:atom:mapping}.

\begin{figure}[htb]
\begin{lstlisting}
    <reaction position="2" id="v2" name="v2" reversible="false">
      <listOfReactants>
        <speciesReference species="A">
          <smtb:carbon position="2" destination="1" occurence="1" species="D"/>
          <smtb:carbon position="1" destination="1" occurence="2" species="D"/>
        </speciesReference>
      </listOfReactants>
      <listOfProducts>
        <speciesReference species="D">
          <smtb:carbon position="1" destination="2" occurence="1" species="A"/>
        </speciesReference>
        <speciesReference species="D">
          <smtb:carbon position="1" destination="1" occurence="1" species="A"/>
        </speciesReference>
      </listOfProducts>
    </reaction>
\end{lstlisting}
\caption{Fragment of the XML file transformed after the first set of XSL transformation: the atom mapping is translated into XML markup
in the namespace \texttt{smtb}}
\label{fig:atom:mapping}
\end{figure}

For each intermediate species, we also add some markup specifying the exhaustive list of its cumomers. For example, for species A in the
branching network, we have the cumomers $\mathrm{A}_{x1}$, $\mathrm{A}_{1x}$, $\mathrm{A}_{11}$ (we don't take into account
$\mathrm{A}_{xx}$ which is equal to 1), and the \verb1<species>1 element corresponding to A is enriched as depicted in Figure
\ref{fig:cumomer:species:list}. Each \verb1<smtb:cumomer>1 element has an id of the form $\mathrm{A}_n$ where $n$ is equal to 
the number whose base 2 representation is equal to the cumomer pattern when replacing the x's by zeros. Each \verb1<smtb:carbon>1 element
denotes a 13 neutrons isotope at position given by the \verb1position1 attribute in the molecule.

\begin{figure}[htb]
\begin{lstlisting}
   <species compartment="default" id="A" name="A" type="intermediate" carbons="2">
      <smtb:cumomer id="A_1" species="A" weight="1" pattern="x1">
        <smtb:carbon position="1"/>
      </smtb:cumomer>
      <smtb:cumomer id="A_2" species="A" weight="1" pattern="1x">
        <smtb:carbon position="2"/>
      </smtb:cumomer>
      <smtb:cumomer id="A_3" species="A" weight="2" pattern="11">
        <smtb:carbon position="1"/>
        <smtb:carbon position="2"/>
      </smtb:cumomer>
    </species>
\end{lstlisting}
\caption{Fragment of the XML file transformed after the first set of XSL transformation: the cumomers of a given intermediate species are
named and their labeling is described using <smtb:carbon> elements.
in the namespace \texttt{smtb}}
\label{fig:cumomer:species:list}
\end{figure}

When we consider the vectors of intermediate species cumomer fractions $\vect{x}_k$ for weights up to 2 for the branching network, we have
$$
\vect{x}_1=\left(\begin{array}{c}\mathrm{A}_{x1},\mathrm{A}_{1x},\mathrm{D}_{1},\mathrm{F}_{x1},\mathrm{F}_{1x}\end{array}\right)^\top,\;
\vect{x}_2=\left(\mathrm{A}_{11},\mathrm{F}_{11}\right)^\top
$$
A redundant enumeration is also generated (see the \verb1<smtb:listOfIntermediateCumomers>1 element on Figure \ref{fig:cumomer:species:globallist}) giving the ordering of all cumomers sorted by weight, allowing to keep the correspondance between
components of vectors $\vect{x}_1$, $\vect{x}_2$ and corresponding species cumomers (this information is needed in the subsequent transformations). A
similar enumeration is also generated for input species cumomers in the \verb1<smtb:listOfInputCumomers>1, giving the correspondance between the
components of vectors $\vect{x}_k^{input}$ and corresponding cumomers :
$$
\vect{x}_1^{input}=\left(\mathrm{A\_{out}}_{x1},\mathrm{A\_out}_{1x}\right)^\top,\;
\vect{x}_2^{input}=\left(\mathrm{A\_{out}}_{11}\right)^\top
$$

\begin{figure}[htb]
\begin{lstlisting}
  <smtb:listOfIntermediateCumomers xmlns:smtb="http://www.utc.fr/sysmetab">
    <smtb:listOfCumomers weight="1">
      <smtb:cumomer id="A_1" species="A" weight="1" pattern="x1" position="1">
        <smtb:carbon position="1"/>
      </smtb:cumomer>
      <smtb:cumomer id="A_2" species="A" weight="1" pattern="1x" position="2">
        <smtb:carbon position="2"/>
      </smtb:cumomer>
      <smtb:cumomer id="D_1" species="D" weight="1" pattern="1" position="3">
        <smtb:carbon position="1"/>
      </smtb:cumomer>
      <smtb:cumomer id="F_1" species="F" weight="1" pattern="x1" position="4">
        <smtb:carbon position="1"/>
      </smtb:cumomer>
      <smtb:cumomer id="F_2" species="F" weight="1" pattern="1x" position="5">
        <smtb:carbon position="2"/>
      </smtb:cumomer>
    </smtb:listOfCumomers>
    <smtb:listOfCumomers weight="2">
      <smtb:cumomer id="A_3" species="A" weight="2" pattern="11" position="1">
        <smtb:carbon position="1"/>
        <smtb:carbon position="2"/>
      </smtb:cumomer>
      <smtb:cumomer id="F_3" species="F" weight="2" pattern="11" position="2">
        <smtb:carbon position="1"/>
        <smtb:carbon position="2"/>
      </smtb:cumomer>
    </smtb:listOfCumomers>
  </smtb:listOfIntermediateCumomers>
  <smtb:listOfInputCumomers xmlns:smtb="http://www.utc.fr/sysmetab">
    <smtb:listOfCumomers weight="1">
      <smtb:cumomer id="A_out_1" species="A_out" weight="1" pattern="x1" position="1">
        <smtb:carbon position="1"/>
      </smtb:cumomer>
      <smtb:cumomer id="A_out_2" species="A_out" weight="1" pattern="1x" position="2">
        <smtb:carbon position="2"/>
      </smtb:cumomer>
    </smtb:listOfCumomers>
    <smtb:listOfCumomers weight="2">
      <smtb:cumomer id="A_out_3" species="A_out" weight="2" pattern="11" position="1">
        <smtb:carbon position="1"/>
        <smtb:carbon position="2"/>
      </smtb:cumomer>
    </smtb:listOfCumomers>
  </smtb:listOfInputCumomers>
\end{lstlisting}
\caption{Fragment of the XML file transformed after the first set of XSL transformation: enumeration of all cumomers sorted by weight and
type (intermediate or input species).}
\label{fig:cumomer:species:globallist}
\end{figure}

\begin{figure}[htb]
\begin{lstlisting}
function [x1,x2,dx1_dv,dx2_dv]=solveCumomers(v,x1_input,x2_input)
  n1=5; // Weight 1 cumomers
  M1_ijv=[1,1,-(v(1)+v(2)+v(3))
  2,2,-(v(1)+v(2)+v(3))
  3,2,v(2)
  3,1,v(2)
  3,3,-(v(4)+v(4))
  4,1,v(1)
  4,2,v(3)
  4,3,v(4)
  4,4,-v(5)
  5,2,v(1)
  5,1,v(3)
  5,3,v(4)
  5,5,-v(5)];
  M1=sparse(M1_ijv(:,1:2),M1_ijv(:,3),[n1,n1]);
  b1_ijv=[1,1,v(6).*x1_input(1,:)
  2,1,v(6).*x1_input(2,:)];
  b1_1=s_full(b1_ijv(:,1:2),b1_ijv(:,3),[n1,1]);

  [M1_handle,M1_rank]=lufact(M1);
  x1=lusolve(M1_handle,-[b1_1]);

  df1_dv_ijv=[1,6,x1_input(1,:)
  1,1,-x1(1,:)
  1,2,-x1(1,:)
  1,3,-x1(1,:)
  2,6,x1_input(2,:)
  2,1,-x1(2,:)
  2,2,-x1(2,:)
  2,3,-x1(2,:)
  3,2,x1(2,:)
  3,2,x1(1,:)
  3,4,-x1(3,:)
  3,4,-x1(3,:)
  4,1,x1(1,:)
  4,3,x1(2,:)
  4,4,x1(3,:)
  4,5,-x1(4,:)
  5,1,x1(2,:)
  5,3,x1(1,:)
  5,4,x1(3,:)
  5,5,-x1(5,:)];
  df1_dv_1=s_full(df1_dv_ijv(:,1:2),df1_dv_ijv(:,3),[n1,6]);
  dx1_dv(:,:,1)=lusolve(M1_handle,-df1_dv_1);

  ludel(M1_handle);

  n2=2; // Weight 2 cumomers
  M2_ijv=[1,1,-(v(1)+v(2)+v(3))
  2,1,v(1)
  2,1,v(3)
  2,2,-v(5)];
  M2=sparse(M2_ijv(:,1:2),M2_ijv(:,3),[n2,n2]);
  b2_ijv=[1,1,v(6).*x2_input(1,:)
  2,1,v(4).*x1(3,:).*x1(3,:)];
  b2_1=s_full(b2_ijv(:,1:2),b2_ijv(:,3),[n2,1]);

  [M2_handle,M2_rank]=lufact(M2);
  x2=lusolve(M2_handle,-[b2_1]);
\end{lstlisting}
\caption{The generated Scilab code for the Branching example (continued on Figure \ref{fig:branchingScilab2})}
\label{fig:branchingScilab1}
\end{figure}

\begin{figure}[htb]
\begin{lstlisting}[firstnumber=61]
  df2_dv_ijv=[1,6,x2_input(1,:)
  1,1,-x2(1,:)
  1,2,-x2(1,:)
  1,3,-x2(1,:)
  2,1,x2(1,:)
  2,3,x2(1,:)
  2,4,x1(3,:).*x1(3,:)
  2,5,-x2(2,:)];
  df2_dv_1=s_full(df2_dv_ijv(:,1:2),df2_dv_ijv(:,3),[n2,6]);
  db2_dx1_ijv=[2,3,x1(3,:).*v(4)
  2,3,x1(3,:).*v(4)];
  db2_dx1_1=sparse(db2_dx1_ijv(:,1:2),db2_dx1_ijv(:,3),[n2,n1]);
  dx2_dv=zeros(n2,6,1);
  dx2_dv(:,:,1)=lusolve(M2_handle,-(df2_dv_1+db2_dx1_1*dx1_dv(:,:,1)));
  ludel(M2_handle);
endfunction

function [cost,grad]=costAndGrad(v)
  [x1,x2,dx1_dv,dx2_dv]=solveCumomers(v,x1_input,x2_input);
  e_label=(C1*x1+C2*x2)-yobs;
  e_flux=E*v-vobs;
  cost=0.5*(sum(delta.*e_flux.^2)+sum(alpha(:,1).*e_label(:,1).^2));
  grad=(delta.*e_flux)'*E+(alpha(:,1).*e_label(:,1))'*(C1*dx1_dv(:,:,1)+C2*dx2_dv(:,:,1));
endfunction

\end{lstlisting}
\caption{The generated Scilab code for the Branching example (begining of code is on Figure \ref{fig:branchingScilab1})}
\label{fig:branchingScilab2}
\end{figure}

\clearpage
\section{Computation of the gradient in the non stationnary case}
The discrete state equation in its cascade form is easily obtained from (\ref{state:discr}) and the definition of $\vect{f}$ as
\begin{equation}
\left(\vect{X}_k-\frac{h}{2}\vect{M}_k\right)\vect{x}_k^{i+1}=
\left(\vect{X}_k+\frac{h}{2}\vect{M}_k\right)\vect{x}_k^i
+\frac{h}{2}
\left(\vect{b}_k(\vect{x}^i)
+\vect{b}_k(\vect{x}^{i+1})
\right),\;1\leq i< N,
\label{state:discr:casc}
\end{equation}
for $k=1\dots n$. We recall that $\vect{b}_k(\vect{x})$ only depends on $\vect{x}_l$ for $l<k$, so that the right-hand side of (\ref{state:discr:casc}) is already known at stage $k$. To obtain $\vect{x}_k^{i+1}$ at each time step $i$ we just have have to solve a sparse linear system with a matrix whose $\vect{LU}$ factors need to be determined only once before the iterations. The cascade structure of discretized state and adjoint state equations is easily recovered. The discretized state equations (\ref{state:discr}) for weights $k=1\dots n$, by

\begin{eqnarray}
\left(\vect{X}_k-\frac{h}{2}\vect{M}_k^\top\right)\vect{p}_k^{N-1}&=&\frac{h}{2}\sum_{l=k+1}^n
\left(\frac{\partial \vect{b}_l}{\partial \vect{x}_k}(\vect{x}^N)\right)^\top \vect{p}_k^{N-1}
 -\frac{\partial I}{\partial \vect{x}_k^N},\\
\left(\vect{X}_k-\frac{h}{2}\vect{M}_k^\top\right)\vect{p}_k^{i-1}&=&
\left(\vect{X}_k+\frac{h}{2}\vect{M}_k^\top\right)\vect{p}_k^i
+\frac{h}{2}\sum_{l=k+1}^n
\left(\frac{\partial \vect{b}_l}{\partial \vect{x}_k}(\vect{x}^i)\right)^\top(\vect{p}_l^{i}+\vect{p}_l^{i-1}),\;1< i< N.
\end{eqnarray}
As in the continous case, the adjoint states $\vect{p}_k$ are obtained in decreasing weight order. The two components of the gradient are finaly obtained by
$$
\left(\frac{dJ(\vect{v},\vect{m})}{d\vect{v}}\right)^\top=-\frac{h}{2}\sum_{i=1}^{N-1}\sum_{k=1}^n\left(\frac{\partial \vect{f}_k}{\partial
\vect{v}}(\vect{x}^{i+1})+\frac{\partial \vect{f}_k}{\partial \vect{v}}(\vect{x}^{i})\right)^\top \vect{p}_k^i,
$$
and
$$
\left(\frac{dJ(\vect{v},\vect{m})}{d\vect{m}}\right)^\top=\sum_{i=1}^{N-1}\sum_{k=1}^n
\left(\frac{\partial}{\partial\vect{m}}\left(\vect{D_k}(\vect{m})\left(\vect{x}_k^{i+1}-\vect{x}_k^{i}\right)\right)\right)^{\top}\vect{p}_k^i,
$$
where for a given vector $\vect{z}\in\mathbb{R}^{n_k}$ the matrix $\frac{\partial}{\partial \vect{m}}(\vect{X}_k(\vect{m})\vect{z})$ is defined by
\begin{equation}
\left(\frac{\partial}{\partial \vect{m}}(\vect{X}_k(\vect{m}) \vect{z})\right)_{ij}=\left\{\begin{array}{rl}
z_i,&\mbox{ if the cumomer fraction }z_i\mbox{ belongs to metabolite }j,\\
0,&\mbox{ otherwise.}
 \end{array}
 \right.
 \label{def:derivm2}
\end{equation}

\end{document}